\documentclass[twoside]{article}
\usepackage{amsmath,amsthm,amssymb}
\usepackage{url}
\usepackage[spanish,english]{babel}
\usepackage{fancyhdr}
\fancypagestyle{plain}{\fancyhf{}%
\fancyhead[L]{\footnotesize{Divulgaciones Matem\'aticas Vol. 16 No. 1%
(2008), pp. \pageref{begin-art}--\pageref{end-art}}}%
}

\fancypagestyle{headings}{\fancyhf{}%
\fancyhead[LE,RO]{\thepage}%
\fancyhead[RE]{Dilcia P\'erez, Yamilet Quintana}%
\fancyhead[LO]{A survey on the Weierstrass approximation theorem}%
\fancyfoot[C]{\footnotesize{Divulgaciones Matem\'aticas Vol. 16 No. 1%
(2008), pp. \pageref{begin-art}--\pageref{end-art}}}}

\newcommand{\CC}{\mathbb{C}}

\newcommand{\NN}{\mathbb{N}}

\newcommand{\pp}{\mathbb{P}}
\newcommand{\RR}{\mathbb{R}}
\newcommand{\GG}{\mathcal{G}}
\newcommand{\ZZ}{\mathbb{Z}}

\newcommand{\es}[1]{\mbox{{\em ess\,lim\/}}#1}

\newcommand{\supp}[1]{\mbox{{\rm supp\/}}#1}
\newcommand{\esup}[1]{\mbox{{\rm ess\,sup\/}}#1}

\newtheorem{teo}{\sc Theorem}[section]

\newtheorem{coro}{\sc Corollary}[section]

\theoremstyle{definition}
\newtheorem{definition}{\sc Definition}[section]
\theoremstyle{remark}
\newtheorem{obs}{\sc Remark}[section]

\begin{document}

\label{begin-art} \pagestyle{headings} \thispagestyle{plain}
\footnote{\hspace{-18.1pt}
Received 2006/07/06. Revised 2007/04/14. Accepted 2007/05/25. 
Research partially supported by  DID-USB under  Grant DI-CB-015-04.}
 \selectlanguage{english}
\begin{center}
{\LARGE\bfseries A survey on the Weierstrass approximation theorem \par }

\vspace{5mm}
{\large\slshape Un estudio sobre  el teorema de aproximaci\'on de
Weierstrass}\\[2mm]

\centerline{{\it Dedicated to the 30th Anniversary of Postgrade in
Mathematics of UCV}}

\vspace{2mm}
{\large Dilcia P\'erez (\url{dperez@ucla.edu.ve})}\\[1mm]
Departamento de Matem\'aticas\\
Universidad Centro Occidental Lisandro Alvarado \\
Barquisimeto-Venezuela.\\[1mm]
{\large Yamilet Quintana (\url{yquintana@usb.ve})}\\[1mm]
Departamento de Matem\'aticas Puras y Aplicadas\\
Universidad Sim\'on Bol\'{\i}var\\
Apartado Postal: 89000, Caracas 1080 A, Venezuela.
\end{center}

\begin{abstract}
The celebrated and famous Weierstrass approximation theorem characterizes the set of  continuous functions on a compact
interval via uniform approximation by algebraic polynomials. This theorem is
the first significant result in Approximation Theory of one real variable and plays a key role in the development of General Approximation Theory. Our aim is to investigate
some new results relative to such theorem, to present a history of the subject, and to introduce some open problems.
\end{abstract}

\selectlanguage{spanish}
\begin{abstract}
El celebrado y famoso teorema de aproximaci\'on de Weierstrass caracteriza al conjunto de las funciones continuas sobre un intervalo compacto v\'{\i}a aproximaci\'on
 uniforme por polinomios algebraicos. Este teorema es el primer resultado significativo en Teor\'{\i}a de Aproximaci\'on de una variable real y juega un rol clave en el desarrollo
de la Teor\'{\i}a de Aproximaci\'on General. Nuestro prop\'osito es investigar algunos de los nuevos resultados relativos a tal teorema,   presentar una historia del tema, e introducir algunos problemas abiertos. \\
{\bf Palabras y frases clave:} Teor\'{\i}a de
Aproximaci\'on, Teorema de Weierstrass, espacios de Sobolev,
espacios de Sobolev con peso, polinomios a valores en $\GG$,
funciones suaves a valores en $\GG$.
\end{abstract}

\selectlanguage{english}
\section{Introduction}
In its most general form, we can say that Approximation
Theory is dedicated to the description of elements  in a topological
space $X$, which can be approximated by elements in a subset $A$ of $X$, that is to say, Approximation Theory allows to characterize the closure of $A$ in $X$.

The first significant results of previous type were those of Karl
Weierstrass (1815-1897),
 who proved in 1885 (when he was 70 years old) the density of algebraic polynomials in the class of continuous
real-valued functions on a compact interval, and the density of
trigonometric polynomials in the class of $2\pi$-periodic
continuous real-valued functions. Such results were -in a sense- a
counterbalance to Weierstrass' famous example of 1861 on the
existence of a continuous  nowhere differentiable function. The
existence of such functions accentuated the need for analytic
rigor in Mathematics, for a further understanding of the nature of
the set of continuous functions, and substantially influenced the
development of analysis. This example represented for some
mathematicians a `lamentable plague' (as Hermite wrote to
Stieltjes on May 20, 1893, see \cite{Bai}), so that,  we can say
the approximation theorems were a panacea. While on the one hand
the set of continuous functions  contains deficient functions, on
the other hand every continuous function can be approximated
arbitrarily well by the ultimate in smooth functions, the
polynomials.

Weierstrass was interested  in Complex Function Theory and in the
ability to represent  functions by power series. The result obtained
in his paper in 1885 should be viewed from that perspective,
moreover, the title of the paper  emphasizes such viewpoint (his
paper was titled {\it On the possibility of giving an analytic
representation to an arbitrary function of real variable}, see
\cite{L}). Weierstrass' perception on analytic functions was of
functions that could be represented by power series.

The paper of Weierstrass (1885) was reprinted in Weierstrass'
Mathematische Werke (collected works)
 with some notable additions, for example a short `introduction'. This reprint appeared in 1903 and
contains the following statement:

{\it The main result of this paper, restricted to the one variable
case, can be summarized as follows:}

{ \it Let $f\in C(\RR)$. Then there exists a sequence $f_1 ,f_2
,\ldots $ of entire functions for which}
$$f(x)=\sum_{i=1}^{\infty}f_{i}(x),$$
{\it for each $x\in\RR$. In addition the convergence of the above sum is
uniform on every finite interval.}

Notice that there isn't mention of the fact that the $f_i$'s may be
assumed to be polynomials.

We state the Weierstrass approximation theorem, not as given in his
paper, but as it is currently stated and understood.

\begin{teo}
\label{t1}
(K. Weierstrass).\\
Given  $f:[a,b]\rightarrow \RR$ continuous and an arbitrary
$\epsilon>0$, there exists an algebraic polynomial $p$ such that
\begin{equation}
\left|f(x)-p(x)\right|\leq \epsilon, \,\, \forall\,\, x\in [a,b].
\end{equation}
\end{teo}

Over the next twenty-five or so years numerous alternative proofs
were given to this result by a roster of some the best analysts of
the period. Two papers of Runge published almost at the same time,
gave another proof of the theorem, but unfortunately, the theorem
was not titled Weierstrass-Runge theorem. The impact of the theorem
of Weierstrass in the world of the Mathematics  was immediate. There
were later proofs of famous mathematicians such as Picard (1891),
Volterra (1897), Lebesgue (1898), Mittag-Leffler (1900), Landau
(1908), de la Valle\'e Poussin (1912). The proofs more commonly
taken at level of undergraduate courses in Mathematics are those of
Fej\'er (1900) and Bernstein (1912) (see, for example \cite{Bar},
\cite{Ch}, \cite{L} or \cite{P}).

We only present a limited sampling of the many results that has been
related with Weierstrass' approximation theorem, which can be found
in Approximation Theory and others areas. We would like to include
all results but the length of this paper would not be suffice
enough. In addition, we don't prove most of the results we quote. We
hope, nonetheless, that the readers will find something here of
interest.

In regards to the outline of the paper, it is as follows:  Section 2
is dedicated to describe  some improvements, generalizations and
ramifications of Weierstrass approximation theorem. Section 3
presents some new results on the subject and to introduce some open
problems.

\section{Improvements, generalizations and ramifications of Weierstrass approximation theorem}

\hspace{.5cm} There are many improvements, generalizations and
ramifications of Weierstrass approximation theorem. For instance, we
could let $f$ be a complex-valued, or take values in a real or
complex vector space of finite dimension; that is easy. We could let
$f$ be a function of several real variables; then it is also easy to
formulate the result. And we could let $f$ be a function of several
complex variables. That would require a more profound study, with
skillful adaptations of both hypothesis and conclusion. Finally, of
course, we may try with functions which take theirs values in an
infinite-dimensional space (see Section 3 below).

Among the results of this kind we can find:

\begin{teo}(Bernstein).\\
Let $f:[0,1]\rightarrow \RR$ bounded, then
$$\lim_{n\to\infty}B_{n}(f;x)=f(x), $$
for each $x\in[0,1]$ where $f$ is continuous. Furthermore, if $f\in
C([0,1])$ then $B_{n}(f;x)$ converges to $f$ uniformly.
\end{teo}

Notice that if $f\in C([a,b])$, we can take $y=\frac{x-a}{b-a}$ for
translate the approximation problem to interval  $[0,1]$ and to use
the Bernstein polynomials $B_{n}(f;y)$.

The theorem of Bernstein besides showing Weierstrass' result, gives
a bonus: an explicit expression for the polynomial approximants to
the function.

 Simultaneously with the study of the approximation for algebraic
 polynomials, also the  approximation  by trigonometrical polynomials was  investigated  -
beginning with  Fourier series-. The problem can generally be
outlined: Given $\{f_{i}\}_{i}$  a sequence of a normed space $X$,
when $\{f_{i}\}_{i}$ is fundamental in $X$?, i.e. when holds that
for every $g\in X$ and  $\epsilon>0$, exist $n$  and a sequence of
scalars  $\{c_{i}\}_{i}$ such that
$$\left\|g-\sum_{i=1}^{n}c_{i}f_{i}\right\|<\epsilon\,\ ?$$

 This question occupied to many mathematicians during the past century. Maybe the oldest and simpler example was proposed by Bernstein in
1912, approximately:  Let $\{\lambda_{i}\}_{i=1}^{\infty}$ a
sequence of distinct positive numbers. When the functions $1,
x^{\lambda_{1}},x^{\lambda_{2}},x^{\lambda_{3}}\ldots $  are
fundamental in  $C([0,1])$?  Also, the answer was surmised by
Bernstein and it is given in the following theorems.

\begin{teo}
\label{teo1}
(M\"untz's first  theorem).\\

Consider the set of functions
$\left\{x^{\lambda_{1}},x^{\lambda_{2}},\ldots\right\}$ where
$-\frac{1}{2}<\lambda_{i}\rightarrow \infty$. It is fundamental in
the least-squares norm on $[0,1]$, if and only if,
$\sum_{\lambda_{i}\not=0}\frac{1}{\lambda_{i}}=\infty$.

\end{teo}

\begin{teo}
\label{teo1a}
(M\"untz's second theorem, 1914).\\
Let $\{\lambda_{i}\}_{i=1}^{\infty}$ a sequence of distinct positive
numbers, such that $1\leq\lambda_{i}\rightarrow \infty$. Then the
set of  functions  $\left\{1,
x^{\lambda_{1}},x^{\lambda_{2}},x^{\lambda_{3}}\ldots\right\}$ is
fundamental in  $C([0,1])$, if and only if,
$\sum_{i=1}^{\infty}\frac{1}{\lambda_{i}}=\infty$.
\end{teo}

By some time this last theorem was called  M\"untz-Szasz theorem,
because Szasz published an independent proof of it in 1916.
Seemingly, Bernstein solved part of the problem. In fact, taking
$\lambda_{i}=i$ we can recover the Weierstrass approximation theorem
for $[0,1] $.

Furthermore, the M\"untz's second theorem is all the more
interesting because it traces a logical connection between two
apparently unrelated facts: the fundamentality of $\left\{1,
x,x^{2},x^{3}\ldots\right\}$ and the divergence of the series of
reciprocal exponents, $\sum_{n=1}^{\infty}\frac{1}{n}$. In fact, if
we wish to delete functions from the set while maintaining its
fundamentality, this divergence is precisely the property that must
be preserved. The reader is referred to \cite{Ch}  for the proofs of
the both theorems.

The combination of  Dual Spaces Theory and Complex Analysis are very
usual in the study of fundamental sequences. Another application of
these techniques is given to solve the problem of Wiener on
traslation of  fundamental sequences, such problem dates from the
thirty's decade of the past century and  is related with the problems
of invariant spaces. This problem has generated almost so much
investigation as M\"untz theorem. The following theorem is a version
of Wiener's theorem.

\begin{teo}
\label{teo2}
(Wiener, approximately 1933).\\
Let $f\in L^{1}(\RR)$. For all $g\in L^{1}(\RR)$ and every $\epsilon
>0$ there are $\{c_{j}\}_{j}, \{\lambda_{j}\}_{j}\subset\RR$, such
that
$$\left\|g(x)-\sum_{j=1}^{n}c_{j}f(x-\lambda_{j})\right\|_{L^{1}(\RR)}<\epsilon,$$
if and only if, the Fourier transform  $f$ doesn't have real zeros,
that is to say,
$$\int_{-\infty}^{\infty}e^{-itx}f(t)dt\not=0, \,\,\forall x\in \RR.$$
\end{teo}

We shall enunciate two  additional theorems  on fundamental
sequences which was published in the 30's of last century and appear
in the Russian edition of the book of Akhiezer \cite{Ach}, and also
in the translation of this book to German, but not in the
translation to English.

\begin{teo}
\label{teo3}
(Akhiezer-Krein. Fundamental sequences of rational functions with simple poles).\\
If $\{a_{j}\}_{j=1}^{\infty}\subset\CC\setminus [-1,1]$, is such
that $a_{j}\not=a_{k}, \,\,\forall j\not=k$, then the sequence
$\left\{\frac{1}{x-a_{j}}\right\}_{j=1}^{\infty}$ is fundamental in
$L^{p}([-1,1])$, $(1\leq p<\infty)$ and in  $ C([-1,1])$, iff
$$\sum_{k=1}^{\infty}\left( 1-\left| a_{k}-\sqrt{a_{k}^{2}-1}\right|\right)=\infty,$$
 where the branch of the function  $f(z)=\sqrt{z}$ considered is the usual, is to say,
$\sqrt{z}>0$ for   $z\in(0,\infty)$.
\end{teo}

\begin{teo}
\label{teo4}
(Paley-Wiener, 1934. Trigonometric series non harmonic).\\
If $\{\sigma_{j}\}_{j=1}^{\infty}\subset\CC$, is such that
$\,\,\sigma_{j}\not=\sigma_{k}, \,\,\forall\,\,\, j\not=k$  and
$|\Im(\sigma_{j})|< \frac{\pi}{2}$ , $\,\,j\geq 1;\,\,$ then the
sequence $\{ e^{ -\frac{\pi}{2} |x| -ix\sigma_{j}
}\}_{j=1}^{\infty}$ is fundamental in  $L^{2}(\RR)$ iff,
$$\sum_{j=1}^{\infty}\frac{\cos( \Im(\sigma_{j}))}{\coth(\Re(\sigma_{j}))}=\infty,$$
where $\Re(z)$, $\Im (z)$ denote, respectively, the real and
imaginary parts of the complex  number $z$.
\end{teo}

Other result which despite its recent discovery has already become
``classical'' and should be a part of the background of every
student of Mathematics is the following.

\begin{teo}
\label{teor41}
(M. H. Stone,  approximately 1947).\\
Let $K$ be a compact subset of $\RR^{n}$ and let $\mathcal{L}$ be a
collection of continuous functions on $K$ to $\RR$ with the
properties:
\begin{enumerate}
\item[i)] If $f, g$ belong to $\mathcal{L}$, then $\sup\{f,g\}$
and $\inf\{f,g\}$ belong to $\mathcal{L}$. \item[ii)] If $a,b\in\RR$
and $x\not= y\in K$, then there exists a function $f$ in
$\mathcal{L}$ such that $f(x)=a$, $f(y)=b$.
\end{enumerate}
Then any continuous function on $K$ to $\RR$ can be uniformly
approximated on $K$ by functions in $\mathcal{L}$.
\end{teo}

The reader is referred to the paper \cite{S1},  in which the
original proof of Stone appears.

Another known result, but not very noted,  is due to E. Hewitt,
\cite{Hew}. This work appeared published in 1947 and it presents
certain generalizations of the Stone-Weierstrass theorem which are
valid in all completely regular spaces.

 There exist another three problems on approximation
which have had a substantial influence on the analysis of last
century:

 \begin{enumerate} \item The problem on polynomial approximation of  Bernstein on whole real line
(between  1912 and 50's). \item To characterize the closure  of
polynomials in families of functions defined on compact sets
$K\subset \CC$ (between  1885 and 50's). It is also necessary to
mention here the  Kakutani-Stone theorem  on the closure of a
lattice of functions real-valued and  Stone-Weierstrass theorem on
the closure of an algebra of functions continuous complex-valued
(see \cite{S1}, \cite{NA1}). \item The Szeg\"o's extremum problem
(between 1920 and 40's).
\end{enumerate}

With respect  the first problem,  it must be treated for non bounded
polynomials in  $\pm\infty$. If let us consider a weight
$w:\RR\rightarrow[0,1]$ (i.e. a non-negative, measurable function)
and we define
$$C_{w}:=\left\{f:\RR\rightarrow\RR, \mbox{ continuous with }\lim_{|x|\to\infty}(fw)(x)=0\right\},$$
with norm
$$\|f\|_{C_{w}}:=\|fw\|_{L^{\infty}(\RR)}.$$

Bernstein wondered when holds that for all $f\in C_{w} $ and every
$\epsilon>0 $,  there exists a polynomial  $p $, such that
\begin{equation}
\label{te1} \|(f-p)w\|_{L^{\infty}(\RR)}<\epsilon\,\,?
\end{equation}

The condition  $fw$ null in  $\pm\infty$ is necessary to give sense
to the problem: we would like $\|pw\|_{L^{\infty}(\RR)}< \infty,\,\,
\forall p\in\pp$, and in particular
$\|x^{n}w\|_{L^{\infty}(\RR)}<\infty,\,\, \forall n\geq 0$. Notice
that  this condition necessarily makes
$\lim_{|x|\to\infty}x^{n}w(x)=0$, then using $(\ref{te1})$ we would
have  $\|fw\|_{L^{\infty}(\RR)}<\epsilon$, therefore $fw$ is null in
$\pm\infty$.

Among those that contributed to the solution of this problem are
Bernstein (1912, 1924),  T. S. Hall (1939, 1950), Dzrbasjan (1947)
and Videnskii (1953). In the case $w$ continuous weight the
solution solution was given by  H. Pollard (1953); other solutions
were given by  Akhiezer (1954), Mergelyan (1956) and Carleson
(1951).

The techniques used in the solution of Bernstein's problem includes
\,Dual \,Spaces Theory, other problems of approximation,
 Entire Functions Theory, analytic  on the half-plane and a lot of
hard work.

 For the approximation problem on compact subset of the complex plane, the situation is as follows:

Let $K$ a compact set of complex plane, we define
$$A(K):=\left\{f:K\rightarrow \CC, \mbox{ continuous in  } K \mbox{ and analytic in
 } K^{\circ}\right\},$$
where $K^{\circ}$ denotes the interior of  $K$ and the norm is given
by
$$\|f\|:=\|f\|_{L^{\infty}(K)}.$$
If we consider
\[ \begin{split}
P(K):= \big\{ f:K\rightarrow \CC,
&\text{ such that there exist polynomials } \{p_n \}  \\
&\text{ with }\lim_{n\to \infty}\|f-p_{n}\|_{L^{\infty}(K)}=0 \big\},
\end{split}
\]
the interesting question is, when $P(K)=A(K)$?

May be this problem was never assigned a name, but in spirit it
should be called the Runge problem (1885),  whose fundamental
contributions to Approximation Theory still can be seen in the course
of Complex Analysis. Among  the interested on the solution of this
problem are J. L. Walsh (1926), Hartogs-Rosenthal (1931),
Lavrentiev (1936) and Keldysh (1945).

\begin{teo}(Mergelyan, 1950's).\\
Let  $K\subset\CC$ compact, the following statements are equivalent:
\begin{enumerate}
\item $P(K)=A(K)$. \item $\CC\setminus K$ is connected.
\end{enumerate}
\end{teo}

With respect to Szeg\"o's extremum problem, it  can be expressed as
follows: Let us consider a non-negative and finite Borel measure
$\mu$ on the unit circle \hbox{$\Pi:=\{z\in\CC: |z|=1\}$} and $p>0$.
The Szeg\"o's extremum problem consists on determining the value
\begin{eqnarray*}\delta_{p}(\mu)&:=&\inf
\left\{\frac{1}{2\pi}\int_{\Pi}\left|P(z)\right|^{p}d\mu(z):P \mbox{
is a polynomial with }
P(0)=1\right\}\\
&=&\inf_{\{c_{j}\}}\left\{\frac{1}{2\pi}\int_{\Pi}\Big|1-\sum_{j\geq
1}c_{j}z^{j}\Big|^{p}d\mu(z) \right\}.
\end{eqnarray*}

Since every polynomial $P(z)$ of degree  $\leq n$ with  $P(0)=1$,
$$Q(z):=z^{n}\overline{P\left(\frac{1}{\overline{z}}\right)}$$
is a monic polynomial of degree $n$,
$\left|Q(z)\right|=\left|P(z)\right|$ whenever $|z|=1$, and  we have
$$\delta_{p}(\mu)= \inf \left\{\frac{1}{2\pi}\int_{\Pi}\left|Q(z)\right|^{p}d\mu(z):Q \mbox{ is monic polynomial}\right\}.$$

This problem, has had more impact in Approximation Theory than
Bernstein's problem   or the Mergelyan  theorem. Its solution and
the techniques developed from it have had a great influence in the
Spaces Hardy Theory and also on the Functions Theory in the 20th.
century.

\begin{teo}
\label{teo5}
(Szeg\"o, 1921).\\
Let $\mu$ absolutely continuous measure with respect to Lebesgue
measure on the unit circle $\Pi$, with
$d\mu(z)=w(e^{i\theta})d\theta$. Then
$$\delta_{p}(\mu)= \exp\left( \frac{1}{2\pi}\int_0^{2\pi}\log (w(e^{i\theta}))d\theta\right).$$
\end{teo}

First Szeg\"o showed this result with $w$  a trigonometric
polynomial and then he used such approximation for extending the
result for $w$ continuous function and also for general case (see
\cite{Sze1}). The generalization to not necessarily absolutely
continuous  measures came later twenty years and is the merit of A.
N. Kolmogorov (1941) for $p=2$, and M. G. Krein (1945) for $p>0$.

\begin{teo}
\label{teo5a}
(Kolmogorov-Krein).\\
If $\mu$  is a non-negative Borel measure on unit circle  $\Pi$,
such that
$$ d\mu(z)=w(z)d\theta+d\mu_{s}(z), \quad z=e^{i\theta}, $$
then
$$\delta_{p}(\mu)= \exp\left( \frac{1}{2\pi}\int_0^{2\pi}\log (w(e^{i\theta}))d\theta\right).$$
\end{teo}

So the singular component $\mu_{s}$ of $\mu$ plays no role at all.

\begin{coro}
\label{cor1}
(Kolmogorov-Krein).\\
Given $p>0$ and a non-negative Borel measure $\mu$  on the unit
circle $\Pi$, such that \hbox{$d\mu(z)=w(z)d\theta+d\mu_{s}(z)$,}
$z=e^{i\theta}$, the following statements are equivalent:
\begin{enumerate}
\item The polynomials are dense in $L^{p}(\mu)$.
\item $\int_0^{2\pi}\log (w(e^{i\theta}))d\theta=-\infty$.
\end{enumerate}
\end{coro}

Therefore, the polynomials are dense in $L^{p}(\mu)$ only in really
exceptional cases. In Orthogonal Polynomials Theory, the condition
$$\int_0^{2\pi}\log (w(e^{i\theta}))d\theta>-\infty$$
is called the Szeg\"o's condition.

In recent years  it has arisen a new focus on the generalizations of
Weierstrass approximation theorem,  which uses the weighted
approximation or approximation by weighted polynomials. More
precisely, if $I$ is a compact interval, the approximation problem
is studied with the norm $L^\infty(I,\, w)$ defined by
\begin{equation}
\|f\|_{L^{\infty}(I,\,w)}:=\esup_{x\in I} |f(x)| w(x)\,,
\end{equation}
where $w$ is a weight, i.e. a non-negative measurable function and
the convention $0\cdot \infty=0$ is used. Notice that (1.1) is not
the usual definition of the $L^\infty$ norm in the context of
measure theory, although it is the correct definition  when we  work
with weights (see e.g. \cite{BO} and \cite{DMS}).

Considering weighted norms $L^\infty(w)$ has been proved to be
interesting mainly because of two reasons: ¯first, it allows to
wider the set of approximable functions (since the functions in
$L^\infty(w)$ can have singularities where the weight tends to
zero); and, second, it is possible to find functions which
approximate f whose qualitative behavior is similar to the one of
f at those points where the weight tends to infinity.

In this case, it is easy to see that $L^{\infty}(I,w)$ and
$L^{\infty}(I)$  are isomorphic, since the map
$\Psi_{w}:L^{\infty}(I,w)\rightarrow L^{\infty}(I)$ given by
$\Psi_{w} (f)=fw$  is a linear and bijective isometry, and
therefore, $\Psi_{w}$ is also homeomorphism, or equivalently,   $
Y\subseteq L^{\infty}(I,w)$,
$\Psi_{w}(\overline{Y})=\overline{\Psi_{w}(Y)}$,  where  we take
each closure with respect  to the norms $L^{\infty}(I,w)$ and
$L^{\infty}(I)$, in each case. Analogously, for all $ A\subseteq
L^{\infty}(I)$,
$\Psi^{-1}_{w}(\overline{A})=\overline{\Psi^{-1}_{w}(A)}$ and
$\Psi^{-1}_{w}=\Psi_{w^{-1}}$. Then  using  Weierstrass' theorem
we have,
\begin{eqnarray}
\label{ecc}
\Psi^{-1}_{w}(\overline{\pp})=\overline{\Psi^{-1}_{w}(\pp)}&=&\{f\in
L^{\infty}(I,\omega): fw\in C(I)\}.
\end{eqnarray}

Unfortunately, in many applications it is not just the isomorphism
type is of interest, since the last equality in (\ref{ecc})
doesn't allow to obtain information on  local behavior of the
functions $f\in L^{\infty}(I,w)$ which can be approximated.
Furthermore,  if \hbox{$f\in L^{\infty}(I,w)$,} then in general
$fw$ is not continuous function, since its continuity depends of
the nature of weight $w$. So, arises the necessary research of the
family of weights for which have sense the approximation problem
by continuous functions. The reader  can find in  \cite{PQRT1} and
\cite{R1} a recent and detailed  study of such problem. It is
precisely in \cite{PQRT1} where appears one of the most general
results known at the present time:

\begin{teo}
\label{teoA} (\cite{PQRT31}, Theorem 2.1).\\ Let $w$ be any weight
and
\[ \begin{split}
H_{0} := \big\{ &f\in L^{\infty}(w): f
\text{ is continuous to the right at every point of } R^{+}(w), \\
 &f \text{ is continuous to the left at every point of  }
R^{-}(w), \\
& \text{for each } a\in S^{+}(w),\quad \es_{x\to a^+} |f(x)-f(a)|\,w(x)=0,\\
& \text{for each } a\in S^{-}(w) ,\quad \es_{x\to a^-} |f(x)-f(a)|\, w(x)=0 \big\}.
\end{split}
\]
Then:
\begin{enumerate}
\item[(a)] The closure of $C(\RR)\cap L^{\infty}(w)$
in $L^{\infty}(w)$ is $H_0$.
\item[(b)] If $w\in L^{\infty}_{loc}(\RR)$, then
the closure of $C^\infty(\RR)\cap L^{\infty}(w)$ in $L^{\infty}(w)$
is also $H_0$.
\item[(c)] If $\supp\,w$ is compact and $w\in L^\infty(\RR)$,
then the closure of the space of polynomials is $H_0$ as well.
\end{enumerate}
\end{teo}

Another special kind of approximation problems arise when we
consider simultaneous approximation which includes to derivatives of
certain functions; this is the case of versions of Weierstrass'
theorem in  weighted Sobolev spaces. In our opinion, very
interesting results  in such  direction appears in \cite{PQRT2}.
Under enough general conditions concerning the vector weight (the so
called type 1) defined in a bounded interval $I$, the authors
characterize to the closure in the weighted Sobolev space $W^{(k,
\infty)} (I ,w)$ of  the spaces of polynomials, $k$-differentiable
functions in the real line, and infinite differentiable functions in
the real line, respectively.

\begin{teo}
\label{teo41}(\cite{PQRT2}, Theorem 4.1).\\
Let us consider a vectorial weight $w=(w_0,\dots,w_k)$ of type $1$
in a compact interval $I=[a,b]$. Then the closure of $\pp \cap
W^{k,\infty}(I,w)$,
$C^\infty(\RR)\cap W^{k,\infty}(I,w)$ and\\
$C^k(\RR) \cap W^{k,\infty}(I,w)$ in $W^{k,\infty}(I,w)$ are,
respectively,

\begin{eqnarray*}
H_1 \!\!\! & :=&  \!\!\! \left\{f\in W^{k,\infty}(I,w):\, f^{(k)}\in
\overline{\pp\cap L^{\infty}(I,w_{k})}^{L^{\infty}(I,w_{k})} \,
\right\} , \\
H_2 \!\!\! & :=&  \!\!\! \left\{f\in W^{k,\infty}(I,w):\, f^{(k)}\in
\overline{C^\infty(I)\cap
L^{\infty}(I,w_{k})}^{L^{\infty}(I,w_{k})} \, \right\} , \\
H_3 \!\!\! & :=&  \!\!\! \left\{f\in W^{k,\infty}(I,w):\, f^{(k)}\in
\overline{C(I)\cap L^{\infty}(I,w_{k})}^{L^{\infty}(I,w_{k})} \,
\right\} .
\end{eqnarray*}
\end{teo}

The reader is referred to \cite{PQRT2} for the proof of this
theorem.

\section{Hilbert-extension of the Weierstrass\\ approximation theorem}
Throughout this  section, let us consider a real and separable
Hilbert space $\GG$, a compact interval $I\subset \RR$,
the space $L^{\infty}_{\GG}(I)$ of all $\GG$-valued essentially
bounded functions, a weakly measurable
function $w:I\rightarrow \GG$, the space $\pp(\GG)$ of all
$\GG$-valued polynomials and the space
$L^{\infty}_{\GG}(I,w)$ of all $\GG$-valued
functions which are bounded whit respect to the norm defined by
\begin{equation}
\label{e2} \|f\|_{L^{\infty}_{\GG}(I,w)}:=\esup_{t\in I}
\|(fw)(t)\|_{\GG},
\end{equation}
where $fw:I\longrightarrow \GG$ is defined as follows: If
$\dim\GG<\infty$, we have the functions $f$ and $w$ can be expressed
as $f=(f_{1},\ldots,f_{n_{0}})$ and $w=(w_{1},\ldots,w_{n_{0}})$,
respectively,  where $f_{k}, \,w_{j}:I\longrightarrow \RR$, for
$j=1,\ldots, n_{0}$, with $n_{0}=\dim \GG$. Then
$$(fw)(t):=
(f_{1}(t)w_{1}(t),\ldots,f_{n_{0}}(t)w_{n_{0}}(t)), \mbox{ for }
t\in I.$$ If $\dim\GG=\infty$, let $\{\tau_{j}\}_{j\in\ZZ_{+}}$ be a
complete orthonormal system, then
$$(fw)(t):=\sum_{j=0}^{\infty}\langle f(t),\tau_{j}\rangle_{\GG}\,\langle
w(t),\tau_{j}\rangle_{\GG}\, \tau_{j}, \mbox{ for } t\in I.$$

It is known that every $\GG$ is a real separable Hilbert space
isomorphic either to $\RR^{n}$ for some $n\in\NN$ or to
$l^{2}(\RR)$. And in both cases $\GG$ has a structure of commutative
Banach algebra, with the coordinatewise operations and identity in
the first case and without identity in the second case (see
\cite{H1}). We use this fact in what follows.

\subsection{ $\GG$-valued functions} Let $\GG$ a real separable Hilbert space. A $\GG$-valued polynomial on $I$ is a function $\phi:I\rightarrow \GG$, such that
$$\phi(t)=\sum_{n\in\NN}\xi_{n}t^{n},$$
where $(\xi)_{n\in\NN}\subset \GG$ has finite support.

Let $\pp(\GG)$ be the space of all $\GG$-valued polynomials on $I$.
It is well known that $\pp(\GG)$ is a subalgebra of the space of all
continuous $\GG$-valued functions on I. We denote by
$L^{\infty}_{\GG}(I)$ the set of all weakly measurable essentially
bounded functions $f:I\rightarrow \GG$. For $1\leq p<\infty$ let
$L^{p}_{\GG}(I)$ denote the set of all weakly measurable functions
$f:I\rightarrow \GG$ such that
$$\int_{I}\|f(t)\|^{p}_{\GG}\,dt<\infty.$$
Then $L^{2}_{\GG}(I)$ is a Hilbert space with respect to inner
product
$$\langle
f,g\rangle_{L^{2}_{\GG}(I)}=\int_{I}\langle
f(t),g(t)\rangle_{\GG}\,dt.$$

$\pp(\GG)$ is also dense in $L^{p}_{\GG}(I)$, for $1\leq p<
\infty$. More details about these spaces can be found in \cite{na}.

In \cite{Q} the author studies the set of functions $\GG$-valued
which can be approximated by $\GG$-valued continuous functions in
the norm $L^\infty_{\GG}(I,w)$, where $I$ is a compact interval,
$\GG$ is a real and separable Hilbert space and $w$ is certain
$\GG$-valued weakly measurable weight, and using the previous
results of \cite{PQRT2} the author obtains a new extension of
Weierstrass' approximation theorem in the context of real separable Hilbert spaces.

\subsection{Approximation in $W_{\GG}^{1,\infty}(I)$}
For dealing with definition of Sobolev space
$W_{\GG}^{1,\infty}(I)$. First of all, notice that we need a
definition of derivative $f^{\prime}$ of a function $f\in
L^{\infty}_{\GG}(I)$, such definition can be introduced by means of
the properties of $\GG$ and the definition of classical Sobolev
spaces as follows.
\newpage
\begin{definition}\mbox{}\\
Given $\GG$ a real separable Hilbert space and $f:I\rightarrow \GG$,
we say that a function $g:I\rightarrow \GG$ is the derivative of $f$
if
$$g\sim\left\{
\begin{array}{l}
(f_{1}^{\prime},\ldots,f_{n_{0}}^{\prime}), \mbox{ if } \dim\GG=n_{0}.\\
\\
\{f_{k}^{\prime}\}_{k\in\NN} \mbox{ with }
\sum_{k\in\NN}|f_{k}^{\prime}(t)|^{2}<\infty, \mbox{ if } \GG \mbox{
is infinite-dimesional. }
\end{array}\right.$$

\end{definition}

So that, we can  define the Sobolev space $W^{1,\infty}_{\GG}(I)$ of
the space of all the  $\GG$-valued essentially bounded functions by
\[ \begin{split}
W^{1,\infty}_{\GG}(I):=\big\{ f\in L^{\infty}_{\GG}(I): &\text{ there exists $f^{\prime}$ in the sense of the Definition 3.1} \\
& \text{ and } f^{\prime}\in L^{\infty}_{\GG}(I)\big\},
\end{split}
\]
with norm
$$\|f\|_{W^{1,\infty}_{\GG}(I)}:=\|f\|_{L^{\infty}_{\GG}(I)}+\|f^{\prime}\|_{L^{\infty}_{\GG}(I)},
$$
\subsection{Interesting problems}

The following problem in this area arises:

If we consider the norm

$$\|f\|_{W^{1,\infty}_{\GG}(I,w)}:=\|f\|_{L^{\infty}_{\GG}(I,w)}+\|f^{\prime}\|_{L^{\infty}_{\GG}(I,w)},
$$
 and we define the weighted Sobolev space of $\GG$-valued functions by means of the relation
$$ f\in W^{1,\infty}_{\GG}(I,w)\Leftrightarrow \left\{
\begin{array}{l} f\in L^{\infty}_{\GG}(I,w),\\
\text{there exists } f^{\prime} \text{ in the sense of the
Definition 3.1}\\
\text{and } f^{\prime}\in L^{\infty}_{\GG}(I,w),\end{array}\right.
$$
then we can ask
\begin{enumerate}
\item What is $ \overline{\pp(\GG)}^{W^{1,\infty}(I,w)}\,?$
\item Which are the most general  conditions on the weight $w$ which allow to characterize $\overline{\pp(\GG)}^{W^{1,\infty}(I,w)}$?
\end{enumerate}
\begin{obs}
Some pieces of information are taken from Internet-based resources
without the URL's.
\end{obs}

\section*{Acknowledgements}
We would like to thank Professors Neptali Romero and Ram\'on Vivas by the offered hospitality while one of the
authors was in UCLA, concluding this paper.

\label{end-art}
\end{document}